\documentclass{amsart}

\usepackage{graphicx}
\usepackage{amssymb}
\usepackage{amsthm}

\input xy
\xyoption{all}

\theoremstyle{plain}
\newtheorem{thm}{Theorem}
\newtheorem{lem}{Lemma}

\newtheorem*{schol}{Scholium}

\theoremstyle{definition}
\newtheorem*{ques}{Question}

\def\T{\mathcal T}
\def\M{\mathcal M}
\def\K{\mathfrak K}
\def\L{\mathfrak L}

\def\h{\mathfrak{H}}

\def\C{\mathbb C}
\def\H{\mathbb H}
\def\R{\mathbb R}

\def\mod{\mathrm{Mod}}
\def\co{\colon\thinspace}

\begin{document}

\title{Totally geodesic boundaries \\ of knot complements}

\author{Richard P. Kent IV}

\address{Department of Mathematics, University of Texas, Austin, TX 78712}
\email{rkent@math.utexas.edu}

\dedicatory{for Kimberly}

\thanks{This work supported in part by a University of Texas Continuing Fellowship.}

\begin{abstract} Given a compact orientable $3$--manifold $M$ whose boundary is a hyperbolic surface and a simple closed curve $C$ in its boundary, every knot in $M$ is homotopic to one whose complement admits a complete hyperbolic structure with totally geodesic boundary in which the geodesic representative of $C$ is as small as you like.  
\end{abstract}

\maketitle
\section{Introduction}

Let $\Sigma$ be a closed hyperbolic surface and let $\h$ be the set of all orientable $3$--manifolds admitting a complete finite volume hyperbolic structure with totally geodesic boundary homeomorphic to $\Sigma$.

By Mostow-Prasad rigidity, there is a map $\widehat \partial$ from $\h$ to the moduli space of $\Sigma$ sending a manifold to the hyperbolic structure that appears on its boundary. It is a theorem of Fujii and Soma that the image of $\widehat \partial$ is dense \cite{fujii}. 

Given a hyperbolic structure $\sigma$, Fujii and Soma construct a sequence of manifolds whose boundaries converge to $\sigma$ using Brooks' theorem \cite{brooks} that the set of hyperbolic structures admitting circle packings is dense. Given $\varepsilon > 0$, a circle packing point $\tau$ is chosen within $\frac \varepsilon 2$ of $\sigma$.  Then a finite volume manifold whose totally geodesic boundary consists of two copies of $\tau$ and two copies of its mirror image is built explicitly from the circle packing.  Two of the boundary components are then identified and hyperbolic Dehn surgeries performed to obtain a manifold with no cusps whose boundary components are within $\frac \varepsilon 4$ of $\tau$ and its mirror image. A number of copies of this manifold are then glued end to end and finally attached to a manifold having totally geodesic boundary merely homeomorphic to $\Sigma$.  A large number may be necessary to ensure that the boundary of this manifold is within $\varepsilon$ of $\sigma$.  In particular, this number tends to infinity as $\varepsilon$ tends to zero. 

With a view towards grasping the relationship between manifold and boundary, we are motivated by the hope that passage to Fujii and Soma's theorem may be made with greater topological control---note, for instance, that as Fujii and Soma's $3$--manifolds approach $\sigma$, their first Betti numbers tend to infinity.

To this end, let $M$ be a compact orientable $3$--manifold whose boundary is homeomorphic to $\Sigma$ and let $\L_M$ be the set of all complements of links in $M$ that lie in $\h$, $\K_M$ the set of all such knot complements.  Fixing a marking  on the boundary of $M$, rigidity again yields a map $\partial$ from $\L_M$ to the Teichm\"uller space of $\Sigma$ sending a manifold to the marked hyperbolic structure that appears on its boundary. 

\begin{ques} Is $\partial \L_M$ dense? Is $\partial \K_M$?
\end{ques}

  We prove the following

\begin{thm}\label{corset} Let $C$ be a simple closed curve in $\partial M$ and let $\varepsilon > 0$. Then every knot in $M$ is homotopic to one whose complement admits a hyperbolic structure with totally geodesic boundary in which the length of the geodesic representative of $C$ is less than $\varepsilon$.
\end{thm}

The proof proceeds as follows. A surface $F$ is chosen in $M$ whose boundary contains $C$.  A theorem of Myers allows us to homotope a given knot to one whose complement admits a hyperbolic structure with totally geodesic boundary. We may chose this knot to intersect $F$ a number of times so that the intersection of $F$ with the resulting knot complement is a quasi-Fuchsian surface.  A pseudo-Anosov mapping class $\psi$ is then found with the property that cutting open along this quasi-Fuchsian surface and precomposing the gluing map with $\psi^n$ yields the complement of a knot in the same homotopy class.  Examination of the final gluing step of the proof of Thurston's Uniformization Theorem and a theorem of J. Brock reveal that these knots are cinching $C$.

\section{Notation and machinery}

\subsection{Topology}

Given two subsets $X$ and $Y$ of a set $Z$ we let $X - Y = X - (X \cap Y)$ denote the set of elements of $X$ that are not elements of $Y$. If $X$ and $Y$ are disjoint and $f\co X \to W$ and $g\co Y \to W$ are two functions, we let $f \cup g$ denote the function $X \cup Y \to W$ whose restrictions to $X$ and $Y$ are $f$ and $g$ respectively. 

If $G$ is a group and $g , h \in G$, we denote the word $g^{-1} h^{-1} g h$ by $[g,h]$.

  Given a topological space $X$, let $|X|$ denote the number of path components of $X$.

  If $M$ is an $n$--manifold and $F$ is a union of transverse properly embedded submanifolds in $M$, we write $M \setminus F = M - \mathrm{int\ nhd}(F)$. If $F$ is a bicollared $(n-1)$--manifold, $\mathrm{nhd} (F) \cong [-1,1] \times F$  and we call $\{-1\} \times F$ and $\{1\} \times F$ the \textit{traces} of $F$ in $M \setminus F$.  For two transverse submanifolds $F$ and $F'$, we often write $F \setminus F'$ for $F \setminus (F \cap F')$.

  A $3$--manifold is said to be \textit{atoroidal} if every one of its embedded incompressible tori is boundary parallel. A properly embedded annulus in a $3$--manifold $M$ is \textit{essential} if it is incompressible and is not isotopic relative to its boundary to an annulus in $\partial M$. A $3$--manifold $M$ is said to be \textit{cylindrical} if it contains an essential annulus. $M$ is said to be \textit{acylindrical} otherwise. 

  A \textit{pared $3$--manifold} is a pair $(M,P)$ where $M$ is an orientable compact irreducible $3$--manifold and $P$ is a subset of $\partial M$ consisting of pairwise nonisotopic (in $\partial M$) embedded incompressible tori and annuli with the property that for $X \in \{S^1 \times S^1, S^1 \times I\}$, every map $(X, \partial X) \to (M, P)$ is homotopic as a map of pairs into $P$. We call $P$ the \textit{pared locus}.

  If $(M, P)$ is a pared manifold with $R = \overline{\partial M -  P} \neq \emptyset$, we let $dM= M \cup_R M$ denote the \textit{double of $M$}---the gluing map is the identity on $R$. A double $dM$ admits an involution $i$ that exchanges the two copies of $M$. We consider $dM$ as a pared manifold with pared locus $dP=P \cup i(P)$.

\subsection{Teichm\"uller space and Kleinian groups}

We refer the reader to \cite{ahl,abi,gard} for more on quasiconformal mappings and Teichm\"uller theory, \cite{mask} for more on Kleinian groups, and \cite{notes} for more on convergence of sequences of hyperbolic manifolds.

  An \textit{analytically finite surface $S'$} is a Riemann surface conformally equivalent to a compact Riemann surface from which a finite set of points has been deleted. The surface $\overline{S'}$ is the Riemann surface obtained from $S'$ by postcomposing every chart with complex conjugation---this canonically reverses the orientation of $S'$. Fix once and for all a compact surface $S$ whose interior is homeomorphic to $S'$. A \textit{marking of $S'$} is a choice of homotopy equivalence $S \to S'$.  The Teichm\"uller space of $S$ parameterizes the marked complete analytically finite hyperbolic structures on the interior of $S$ up to isotopy and is denoted $\T(S)$.  

  Let $S_1$ and $S_2$ be analytically finite surfaces equipped with markings $f_1\co S \to S_1$ and $f_2\co S \to S_2$. A homeomorphism $f\co S_1 \to S_2$ is \textit{$K$--quasiconformal} if it is absolutely continuous on lines and $|f_{\bar z}| \leq k|f_z|$ in every local coordinate $z$, where $k = (K-1)/(K+1) < 1$. The minimum possible value of $K$ for which $f$ is $K$--quasiconformal is called the \textit{dilatation of $f$}, and is denoted $K[f]$.  The \textit{Teichm\"uller distance} between $S_1$ and $S_2$ is then
\[
d_\T(S_1,S_2) =  \frac{1}{2} \inf \log K[f]
\]
where the infimum is taken over all quasiconformal maps $f\co S_1 \to S_2$ homotopic to $f_2 \circ f_1^{-1}$. There is a unique extremal quasiconformal map $S_1 \to S_2$  realizing the above distance, called the \textit{Teichm\"uller mapping}. 

  We follow the convention that if $S$ is disconnected, then its Teichm\"uller space is the product of the Teichm\"uller spaces of its components equipped with the metric obtained by taking the maximum of the Teichm\"uller metrics on the factors. If $S$ is of genus $g$ and connected, $\T(S)$ is homeomorphic to an open Euclidean ball of dimension $6g-6 + 2|\partial S|$. 

The mapping class group $\mod(S) = \pi_0 (\mathrm{Homeo}^+(\mathrm{int}\ S))$ acts properly discontinuously on $\T(S)$ by pulling back hyperbolic structures.  The moduli space $\M(S)$ is the quotient of $\T(S)$ by this action. A mapping class is said to be \textit{reducible} if it fixes a collection $\mathcal C$ of disjoint nonperipheral simple closed curves up to isotopy. Elements of $\mathcal C$ are called \textit{reducing curves}. A mapping class $\varphi$ is said to be \textit{pseudo-Anosov} if no nontrivial power of $\varphi$ is reducible. 

We let $\mathcal P(S)$ denote the mapping class group of isotopy classes of homeomorphisms of $S$ that do not permute the components of $\partial S$. Since we are not requiring boundary components to be fixed pointwise, we may regard $\mathcal P(S)$ as a subgroup of $\mod(S)$.

Let $\H^3$ denote hyperbolic $3$--space, $S^2_\infty = \widehat \C$ its sphere at infinity. 

  For our purposes, a \textit{Kleinian group $\Gamma$} is a discrete subgroup of $\mathrm{PSL}_2(\C) \cong \mathrm{Isom}^+(\H^3)$.  We write $M_\Gamma= \H^3 /\Gamma$. The \textit{limit set $\Lambda_\Gamma$ of $\Gamma$} is the smallest nonempty closed subset of $S^2_\infty$ that is invariant under the action of $\Gamma$. 

Let $H_\Gamma$ denote the hyperbolic convex hull of $\Lambda_\Gamma$ in $\H^3$.  A Kleinian group $\Gamma$ is said to be \textit{geometrically finite} if there is an open neighborhood of $H_\Gamma/\Gamma$ in $M_\Gamma$ of finite volume. 

  The \textit{domain of discontinuity of a Kleinian group $\Gamma$} is the set $\Omega_\Gamma = S^2_\infty - \Lambda_\Gamma$, on which $\Gamma$ acts properly discontinuously.  By the Ahlfors Finiteness Theorem, the surface $\Omega_\Gamma/\Gamma$ is analytically finite whenever $\Gamma$ is nonelementary---see \cite{kapo}, Section 4.19, and the references cited there. This surface inherits a conformal structure from that of $\widehat \C$. By the Uniformization Theorem, $\Omega_\Gamma/\Gamma$ also comes equipped with a hyperbolic structure.  The manifold $\dot M_\Gamma = (H^3 \cup \Omega_\Gamma)/\Gamma$ is the \textit{Kleinian manifold associated to $\Gamma$}.   

  A \textit{hyperbolic structure on a pared manifold $(M, P)$} is a discrete faithful representation $\pi_1(M,*) \to \mathrm{PSL}_2(\C)$ with image $\Gamma$  in which the images of the subgroups corresponding to the components of $P$ consist entirely of parabolic elements and such that $\H^3/\Gamma \cong \mathrm{int}\ M$. The Kleinian manifold $\dot M_\Gamma$ is homeomorphic to $M - P$.  A hyperbolic structure is \textit{geometrically finite} if $\Gamma$ is geometrically finite.  

A hyperbolic structure is said to be a \textit{hyperbolic structure with totally geodesic boundary} if the images of the subgroups corresponding to the components of $\partial M - P$ are Fuchsian.

Two hyperbolic structures are deemed \textit{equivalent} if their images are conjugate in $\mathrm{PSL}_2(\C)$, and we let $\mathrm{H}(M,P)$ denote the set of equivalence classes.  We continue to refer to equivalence classes as hyperbolic structures.   

  Let $\mathrm{AH}(M,P)$ denote the set $\mathrm{H}(M,P)$ equipped with the topology induced by the inclusion 
\[
\mathrm{H}(M,P) \subset \mathrm{Hom}(\pi_1(M,*), \mathrm{PSL}_2 (\C))/\mathrm{PSL}_2(\C)
\]
 equipped with the quotient of the compact--open topology---this is the topology induced by convergence on a set of generators after conjugating and is called the \textit{algebraic topology}.

  Let $M$ and $N$ be two manifolds equipped with hyperbolic metrics $d_M$ and $d_N$ respectively and let $K > 1$. A diffeomorphism $f\co M \to N$ satisfying
\[
K^{-1} d_M(x,y) \leq d_N(f(x),f(y)) \leq K d_M(x,y),
\]
for all $x$ and $y$ in $M$ is a called a \textit{$K$--quasi-isometry}.

\subsection{Quasi-Fuchsian groups and iteration near a Bers slice}

A Kleinian group $\Gamma$ is \textit{quasi-Fuchsian} if $\Lambda_\Gamma$ is a topological circle and $\Gamma$ preserves each component $\Omega_1$ and $\Omega_2$ of $\Omega_\Gamma$ setwise. The surfaces $S_1 = \Omega_1/\Gamma$ and $S_2 = \Omega_2/\Gamma$ are homeomorphic by an orientation-reversing homeomorphism and each comes equipped with a hyperbolic structure of its own.  So, associated to $\Gamma$ is a point in $\T(S) \times \T(\overline S)$ once a marking $S \to S_1$ has been chosen.  We call $S$ the \textit{underlying surface of $\Gamma$}. The manifold $\H^3 / \Gamma $  is called the \textit{quasi-Fuchsian manifold associated to $\Gamma$} and is homeomorphic to $S \times \R$. Let $\mathrm{QF}(S) \subset \mathrm{AH}(S \times \R)$ denote the space of quasi-Fuchsian groups. The Ahlfors-Bers Simultaneous Uniformization Theorem \cite{simulbers} provides a homeomorphism $Q\co \T(S) \times \T(\overline S) \to \mathrm{QF}(S)$, natural in the sense that $Q^{-1}(\Gamma)$ is the pair of surfaces $\Omega_\Gamma/\Gamma$. If $A \subset \T(S)$, and $B \subset \T(\overline S)$, we write $Q(A,B) = Q(A \times B)$ and we abbreviate $Q(\{a\},\{b\})$ by writing $Q(a,b)$.

We will need the following

\begin{lem}\label{close} Let $A \subset \T(S)$ and $B \subset \T(\overline{S})$ be compact sets.  Then there is a $k$ such that for each $\varphi \in \mod(S)$, all quasi-Fuchsian manifolds associated to elements of $Q(\varphi A,B)$ are $k$--quasi-isometric. 
\end{lem}
\begin{proof} Since $A$ is compact and $\mod(S)$ acts on $\T(S)$ by isometries of the Teich-m\"uller metric, the diameter of $\varphi A$ in this metric is bounded independent of $\varphi$.

  Fix $b \in B$, $a_1, a_2 \in A$. Let $\Omega_i \subset \Omega_{Q(\varphi a_i,b)}$ be the components  corresponding to $\varphi a_i$, $\widehat \Omega_i \subset \Omega_{Q(\varphi a_i,b)}$ the components corresponding to $b$ and let $h_i\co U \to \Omega_i$ and $g_i\co L \to \widehat \Omega_i$ be Riemann mappings from the upper and lower halfplanes in $\widehat \C$, respectively, for $i \in\{1,2\}$. Let $\Gamma_i$ be the Fuchsian group uniformizing $\varphi a_i$ for each $i$. We may lift the Teichm\"uller mapping $\varphi a_1 \to \varphi a_2$ to obtain a quasiconformal mapping $\omega \co U \to U$ such that $\omega^{-1} \Gamma_2 \omega = \Gamma_1$ and such that $h_1 \omega^{-1} h_2^{-1} \cup g_1 g_2^{-1}$ extends to a quasiconformal mapping $w\co \widehat \C \to \widehat \C$ with $w Q(\varphi a_1,b) w^{-1} = Q(\varphi a_2, b)$---compare Section 6 of \cite{bers}. 

  Since the diameter of $\varphi A$ is uniformly bounded, the dilatation of this map is also uniformly bounded as the $a_i$ range over $\varphi A$.  

  Symmetry allows us to vary $b$ in $B$, and so there is a constant $k'$ depending only on $A$ and $B$ such that any $Q(\varphi a,b)$ and $Q(\varphi a',b')$ are $k'$--quasiconformally conjugate Kleinian groups whenever $a,a' \in A$ and $b,b' \in B$. Such quasiconformal conjugations extend to equivariant $k$--quasi-isometries between universal covers of quasi-Fuchsian manifolds---see Theorem 2.5 of \cite{mc1}---with $k$ uniform over $Q(\varphi A, B)$ and independent of $\varphi$. By equivariance, these $k$--quasi-isometries descend to $k$--quasi-isometries between quasi-Fuchsian manifolds.
\end{proof}

   Let $\varphi \in  \mod(S)$ be reducible of infinite order.  By passing to a suitable power, we may assume that $\varphi$ fixes every reducing curve so that $\varphi$ is pseudo-Anosov or trivial when restricted to any component of the complement of the reducing curves and we do so.  Consider the sequence of quasi-Fuchsian groups $Q(\varphi^n a,b)$ and their quasi-Fuchsian manifolds $M_n$.  It is a theorem of J. Brock \cite{brock} that the lengths of the reducing curves in $M_n$ tend to zero as $n$ approaches infinity. 

  Lemma \ref{close} will allow us to apply Brock's theorem to sequences of quasi-Fuchsian groups whose $n^\mathrm{th}$ terms lie in sets $Q(\varphi^n A, B)$.

\subsection{Uniformization and the skinning map}

We review here the final gluing step of the proof of Thurston's Uniformization Theorem for Haken Manifolds and refer the reader to \cite{morg}, \cite{kapo} and \cite{mask} for more information regarding this theorem, the Bounded Image Theorem, and Maskit combination.

Let $S$ be a two-sided properly embedded surface in a pared 3--manifold $(M,P)$. Let $M'=M\setminus S$ and let $S_1$ and $S_2$ denote the traces of $S$ in $M'$.  Let $f\co S_1 \to S_2$ be the gluing map such that $M = M'/f$. We mark $S_1$ with $f_1 \co S \to S_1$ and $S_2$ with $f \circ f_1$. Precomposing $f$ with a representative of a mapping class $\varphi \in \mod(S_1)$ yields a new gluing map $g$ and we let $M(\varphi) = M'/g$.

  Given a geometrically finite hyperbolic structure $\rho$ on a pared manifold $(M,P)$ with image $\Gamma$, we obtain a point in $\T(\partial M-P)$ by considering $\Omega_\Gamma/\Gamma$. Conversely, by a theorem of Bers \cite{quasibers}, if $(M,P)$ admits such a hyperbolic structure, then the set of all geometrically finite hyperbolic structures on $(M,P)$ is parameterized by $\T(\partial M - P)$.

  So, let $(M,\partial M)$ be a pared  $3$--manifold that does not admit a Seifert fibration and suppose that $S$ is a two-sided properly embedded incompressible surface such that $(N,P)= (M\setminus S, \partial M \setminus \partial S)$ is acylindrical and admits no Seifert fibration. Let $S_1$ and $S_2$ denote the traces of $S$ in $N$ equipped with markings $f_i \co S \to S_i$ as above. Then $(N,P)$ admits a geometrically finite hyperbolic structure $\rho$ by the Uniformization Theorem for Haken Manifolds and so the set of such structures is parameterized by $\T(\partial N - P) \cong \T(S_1) \times \T(S_2)$.  

  Given an element of $\T(S_1) \times \T(S_2)$, consider the geometrically finite structure on $(N,P)$ to which it is associated. Now, for any pared $3$--manifold $(M',P')$ equipped with a geometrically finite hyperbolic structure, the subgroup corresponding to a component of $\partial M' - P'$ is quasi-Fuchsian or else $M'$ is cylindrical, by Proposition 7.2 of \cite{morg} and the Annulus Theorem. Since $(N,P)$ is acylindrical, the images $\Gamma_1 = \rho(\pi_1(S_1))$ and $\Gamma_2=\rho(\pi_1(S_2))$ are quasi-Fuchsian groups.  Now,  $\Gamma_1$ determines a point $(a, \sigma_1(a)) \in \T(S_1) \times \T(\overline{S_1})$.  The first coordinate is the conformal structure that appears on the $S_1$ boundary component of $(N,P)$ and $\sigma_1(a)$ is called the \textit{conformal structure that $S_1$ hides}. The map
\[
\sigma_1 \co \T(S_1) \longrightarrow \T(\overline{S_1})
\]
is called the \textit{skinning map (associated to $N$)}.  

  Since $(N,P)$ is acylindrical, $\sigma_1$ has precompact image---this is the content of Thurston's Bounded Image Theorem. We have a similar discussion and a map $\sigma_2$ associated to $S_2$. For $i\in \{1,2\}$, let $\overline{\sigma}_i \co \T(S_i) \to \T(S_i)$ denote $\sigma_i$ followed by the map induced by reversing the orientation of $\overline{S_i}$. 

  Let $\varphi \in \mod(S_1)$. Precomposing the gluing map $f$ by a representative of $\varphi$ yields two sequences of maps:
\[
\xymatrix{
\T(S_1) \ar[r]^{\overline{\sigma}_1} & \T(S_1) \ar[r]^{\varphi_*} & \T(S_1) \ar[r]^{f_*} &\T(S_2)
}
\]
and
\[
\xymatrix{
\T(S_2) \ar[r]^{\overline{\sigma}_2} & \T(S_2) \ar[r]^{f_*^{-1}} & \T(S_1) \ar[r]^{\varphi_*^{-1}} &\T(S_1)
}
\]
To obtain a hyperbolic structure for $(M,\partial M)$, an application of one of two of the Maskit Combination Theorems is applied.  Satisfaction of the hypotheses is easily translated into the existence of a point $(x,y) \in \T(S_1) \times \T(S_2)$ such that
\[
\overline{\sigma}_2(y) = f_* \varphi_* (x)
\]
and
\[
\overline{\sigma}_1(x) = \varphi_*^{-1} f_*^{-1}(y)
\]
We call this point the \textit{solution to the gluing problem for $f \varphi$}. The existence of such a point means that there is a geometrically finite hyperbolic structure on $(N,P)$ so that the conformal structure at infinity associated to $S_1$ is carried by our gluing map to the conformal structure that $S_2$ hides, and vice versa.

\section{Knots}\label{laces}

Let $M$ be a compact orientable $3$--manifold whose boundary is homeomorphic to $\Sigma$ and let $C$ be a simple closed curve in $\partial M$. If $C$ is separating in $\partial M$, let $F'$ be a surface in $\partial M$ that it bounds.  If $C$ is nonseparating, it is a boundary component of a pair of pants $F'$ whose boundary components are pairwise nonisotopic in $\partial M$. In either case, let $F$ be a properly embedded surface in $M$ isotopic to $F'$.

\begin{lem}\label{rope}  Every homotopy class of simple closed curves in $M$ contains a knot $K$ such that $M \setminus K$ admits a hyperbolic structure with totally geodesic boundary and $M \setminus(F \cup K)$ is boundary incompressible and acylindrical.
\end{lem}
\begin{proof} Let $\gamma$ be a homotopy class and let $k$ be a simple closed curve representing $\gamma$ that intersects $F$ transversely in at least two points.  Consider $k'= k \setminus F \subset M \setminus F$.  By Theorem 1.1 of \cite{myers}, there is a $1$--manifold $k''$ in the homotopy class of $k'$ relative to its boundary such that each component of $(M \setminus F) \setminus k''$ is boundary incompressible, acylindrical and atoroidal. Let $K$ be the knot in $M$ corresponding to $k''$.

The knot exterior $M \setminus K$ is irreducible, atoroidal and acylindrical by Lemma 2.1 of \cite{myers}.  The double $d(M \setminus K)$ admits a hyperbolic structure by the Uniformization Theorem for Haken Manifolds and we equip $d(M - K)$ with the associated hyperbolic metric. The natural orientation-reversing involution on $d(M - K)$ is homotopic to an isometry that preserves a totally geodesic surface, by Mostow-Prasad rigidity.  A cut and paste argument demonstrates that this surface is isotopic to $\partial (M-K)$ and cutting the double open again yields a hyperbolic structure with totally geodesic boundary on $M \setminus K$. 
\end{proof}

\begin{schol}\label{stitches} Given $M$, $F$ and $K$ as in Lemma \ref{rope}, cutting $M \setminus K$ along $F_K=F \setminus K$ and regluing the traces of $F_K$ in a different fashion yields a hyperbolic manifold with totally geodesic boundary no matter what mapping class we choose for the gluing.
\end{schol}
\begin{proof} The conclusion of Lemma 2.1 of \cite{myers} holds for any gluing map.
\end{proof}

\section{Proof of Theorem \ref{corset}}

\begin{proof}[Proof of Theorem \ref{corset}] Let $M$, $F$ and $K$ be as in Lemma \ref{rope} and let $F_K=F\setminus K$. Let $(X_K,P) = (M \setminus K, \partial X_K - \partial M)$, $Y_K = dX_K$, $S = dF_K \subset Y_K$, $Z_K = Y_K \setminus S$. Let $S_1$, $S_2$ denote the traces of $S$ in $Z_K$, $f \co S_1 \to S_2$ the gluing map such that $Z_K/f = Y_K$. We mark $S_1$ with $f_1 \co S \to S_1$, $S_2$ with $f \circ f_1$. Note that $Z_K$ is acylindrical as no component of $\partial M \setminus \partial F$ is an annulus. 

Let $\{\delta_1, \ldots,\delta_{\ell+1}\}$ be the set of boundary components of $F_K$ and let $\widehat F_K$ denote the surface obtained from $F_K$ by capping off a boundary component with a disk.  Let $\{\gamma_1, \ldots, \gamma_\ell\}$ be the set of boundary components of $\widehat F_K$. For each $\delta_i$, there is a short exact sequence
\[
1 \to \pi_1(\widehat F_K) \to \mathcal P(F_K) \stackrel{\phi_i}{\longrightarrow} \mathcal P(\widehat F_K) \to 1
\]
where an element of $\pi_1(\widehat F_K)$ corresponds to the mapping class that ``spins''  the associated boundary component about the homotopy class and $\phi_i$ ``forgets'' $\delta_i$, see \cite{birman2}.  Note that for any $\psi \in \bigcap_{i=1}^{\ell+1}  \mathrm{ker}\thinspace \phi_i$, $X_K(\psi)$ is the complement of a knot in the same homotopy class as $K$. We claim that $\bigcap_{i=1}^{\ell+1}  \mathrm{ker}\thinspace \phi_i$ is nontrivial. To see this, choose a basis $\{x_1, \ldots, x_m\}$ for the free group $\pi_1(\widehat F_K)$ so that a representative for $x_i$ is freely homotopic to $\gamma_i$ whenever $1\leq i \leq \ell-1$ and a representative for 
\[
v= x_1 \ldots x_{\ell-1}[x_\ell,x_{\ell+1}][x_{\ell+2},x_{\ell+3}] \ldots [x_{m-1},x_m]
\]
 is freely homotopic to $\gamma_\ell$---such a basis is easily obtained from the standard handle decomposition of $\widehat F_K$.  Since $\pi_1(\widehat F_K)$ is nonabelian, the mapping class associated to the commutator 
\[
w = \big [ [\ldots [ [x_1,x_2],x_3], \ldots ,x_m], v \big ]
\]
 is a nontrivial element of $\bigcap_{i=1}^{\ell+1}  \mathrm{ker}\thinspace \phi_i$. 

We claim that $w$ represents a pseudo-Anosov mapping class $\psi$. Suppose that this is not the case. By theorem 2 of \cite{kra}, $w$ has a representative on $\widehat F_K$ that does not fill $\widehat F_K$.  This implies that $w$ is conjugate into a proper free factor of $\pi_1(\widehat F_K)$.  It is easily verified, using techniques due to J. H. C. Whitehead \cite{whitehead}, that this is not the case, as the Whitehead graph associated to $w$ has no cut vertex, see \cite{stallings}.

  Let $\varphi = d \psi$ be the double of $\psi$.  

  We consider $X_K(\psi^n)$, $Y_K(\varphi^n)$, and $Z_K$ as pared manifolds with pared loci  $\partial X_K(\psi^n) - \partial M$, $\partial Y_K(\varphi^n)$, and $\partial Z_K \setminus (S_1 \cup S_2)$ respectively, and henceforth suppress the paring in the notation. 

 The manifold $X_K(\psi^n)$ admits a hyperbolic structure with totally geodesic boundary by the Scholium and $Y_K(\varphi^n)$ admits a hyperbolic structure by doubling.  

  The latter structure may be obtained from the unique solution to the gluing problem for $f \varphi^n$ and so we turn our attention to the double $Y_K$ cut open, $Z_K$.  

  We have the sequences obtained by precomposing $f$ by $\varphi^n$: 
\[
\xymatrix{
\T(S_1) \ar[r]^{\overline{\sigma}_1} & \T(S_1) \ar[r]^{\varphi_*^n} & \T(S_1) \ar[r]^{f_*} &\T(S_2)
}
\]
and
\[
\xymatrix{
\T(S_2) \ar[r]^{\overline{\sigma}_2} & \T(S_2) \ar[r]^{f_*^{-1}} & \T(S_1) \ar[r]^{\varphi_*^{-n}} &\T(S_1)
}
\]
Since $Z_K$ is acylindrical, the Bounded Image Theorem demands that the images of the $\overline{\sigma}_i$ be contained in compact $A_i$ for each $i \in \{1,2\}$. The solution to the gluing problem for $f$ lies in the compact set 
\[
f^{-1}_*A_2 \times f_* A_1 \subset \T(S_1) \times \T(S_2)
\]
 and the solution to the gluing problem for $f \varphi^n$ lies in the compact 
\[
\varphi_*^{-n} f_*^{-1} A_2 \times f_*\varphi_*^n A_1 \subset \T(S_1) \times \T(S_2)
\]
Let $A=f_*^{-1} A_2$ and let $B \supset \sigma_1(\T(S_1))$ be compact. Write $\Phi = \varphi^{-1}$. In the resulting hyperbolic structure on $Y_K(\varphi^n)$,  $\pi_1(S)$ is a quasi-Fuchsian group $Q(\Phi^n a_n, b_n)$ trapped in $Q(\Phi^n  A, B)$. 

  Let $a \in A$, $b \in B$.   Since $\Phi$ is reducible and pseudo-Anosov when restricted to either component of $S\setminus \partial F$, the lengths of the geodesic representatives of $\partial F$ in the quasi-Fuchsian manifolds associated to the $Q(\Phi^n a, b)$ tend to zero as $n$ grows, by Theorem 4.5 of \cite{brock}. By Lemma \ref{close}, the quasi-Fuchsian manifolds associated to $Q(\Phi^n a_n, b_n)$ and $Q(\Phi^n a, b)$ are quasi-isometric independent of $n$, and so the lengths of the geodesic representatives of $\partial F$ in the former  manifolds also tend to zero as $n$ grows. In particular, the lengths of the geodesic representatives of $C$ in the $Y_K(\varphi^n)$ tend to zero. Now $\partial X_K(\psi^n)$ is totally geodesic in $Y_K(\varphi^n)$ and so $C$ is shrinking in the $\partial X_K(\psi^n)$ as well.
\end{proof}

In the case that $C$ is nonseparating, we have cinched $C$ at the expense of cinching two other curves in $\partial M$.  This may be remedied at the expense of the homotopy class of the knot.  Consider two disjoint parallel pushoffs $C'$ and $C''$ of $C$ just inside $\partial M$.  Pick an annulus $F$ between $C$ and $C'$ disjoint from $C''$ and proceed as above to obtain a knot $K$ in $M - (C' \cup C'')$ in whose complement $C$ is short.  Performing $\frac 1 m$--filling on the cusp corresponding to $C'$ and $- \frac 1 m$--filling on the cusp corresponding to $C''$ has the effect of shearing $K$ along an annulus whose boundary is $C' \cup C''$, yielding a knot in $M$.  As $m$ tends to infinity, these knot complements converge geometrically to $M - (K \cup C' \cup C'')$, see Proposition E.6.29 of \cite{bp}, and so for large $m$, $C$ is short in these manifolds.  Note that this construction is easily adapted to cinch a collection of disjoint simple closed curves.

\section{Geometric flexibility and isolation}

In \cite{neureid}, Neumann and Reid construct infinitely many examples of one-cusped hyperbolic $3$--manifolds that contain totally geodesic surfaces that are \textit{geometrically isolated} from the cusp in the sense that performing any hyperbolic Dehn filling leaves the hyperbolic structure on the totally geodesic surface unchanged. In contrast, Fujii \cite{fuj} constructs an explicit one-cusped hyperbolic manifold $M$ with totally geodesic boundary with the property that the manifolds obtained from $M$ by filling exhibit infinitely many hyperbolic structures on their boundaries and in \cite{fuj2}, one-cusped manifolds with totally geodesic boundaries are constructed for which the map from the hyperbolic Dehn filling space to the Teichm\"uller space of the boundary is actually an embedding near infinity. 

We have the

\begin{thm} Let $M$ be a $3$--manifold that admits a hyperbolic structure with totally geodesic boundary homeomorphic to $\Sigma$.   Then $M$ contains a hyperbolic knot $K$ such that $\partial M$ and $\partial (M-K)$ are as far apart as we like in $\M(\Sigma)$.
\end{thm}
\begin{proof} Theorem \ref{corset} provides a knot $K$ such that the injectivity radius of $\partial (M-K)$ is as small as we like, so we may choose $K$ so that $\partial (M-K)$ is as far as we like from $\partial M$ in $\M(\Sigma)$---since the set of surfaces in $\M$ with injectivity radius greater than or equal to $\varepsilon >0$ is compact \cite{mum}. 
\end{proof}

A hyperbolic manifold obtained by performing Dehn filling on a hyperbolic manifold is said to be obtained by \textit{hyperbolic Dehn filling} if the core of the filling torus is a geodesic---see section E.4 of \cite{bp}. Unlike the fillings considered in \cite{neureid}, \cite{fuj}, and \cite{fuj2}, the filling on $M - K$ yielding $M$ is typically not a hyperbolic Dehn filling---for example, we may take $K$ to be null-homotopic in $M$. So, though $M - K$ admits a filling that moves the boundary, it could be that every hyperbolic Dehn filling leaves the hyperbolic structure on the boundary alone.

\section*{Acknowledgements}
 The author thanks his advisor Cameron Gordon, Jeff Brock, Chris Leininger, and Alan Reid for useful conversations. He thanks the referee for thoughtful comments.


\begin{thebibliography}{99}

\bibitem{abi} W. Abikoff, The real analytic theory of Teichm\"uller space, Lecture Notes in Math, vol. 820, Springer, 1980.

\bibitem{ahl} L. V. Ahlfors, Lectures on quasiconformal mappings, Van Nostrand, 1966.

\bibitem{bp} R. Benedetti and C. Petronio, Lectures on hyperbolic geometry, Springer-Verlag 1992.

\bibitem{bers} L. Bers, On boundaries of Teichm\"uller spaces and on Kleinian groups: I.  Ann. of Math., 91 (1970), 570--600.

\bibitem{quasibers} $\underline{\ \ \ \ \ \ \  \ \ \ \ }\thinspace$, Quasiconformal mappings and Teichm\"uller's theorem, Analytic Functions (R. Nevanlinna et al., eds), 89--119. Princeton UP 1960.

\bibitem{simulbers} $\underline{\ \ \ \ \ \ \  \ \ \ \ }\thinspace$, Simultaneous uniformization, Bull. AMS 66 1960 94--97.

\bibitem{birman2} J. Birman, Mapping class groups and their relationship to braid groups, Comm. Pure and Applied Math. 22(1969) 213-238.

\bibitem{brock} J. F. Brock, Iteration of mapping classes and limits of hyperbolic 3--manifolds, Invent. Math. 143, 523--570 (2001).

\bibitem{brooks} R. Brooks, Circle packings and co-compact extensions of Kleinian groups, Invent. Math. 86, 461--469 (1986).

\bibitem{notes} R. D. Canary, D. B. A. Epstein and P. Green, Notes on notes of Thurston. Analytical and geometric aspects of hyperbolic space (Coventry/Durham, 1984), 3--92, London Math. Soc. Lecture Note Ser., 111, Cambridge UP, 1987. 

\bibitem{fuj} M. Fujii, On totally geodesic boundaries of hyperbolic $3$--manifolds, Kodai Math. J. 15 (1992), 244--257

\bibitem{fuj2} M. Fujii and S. Kojima, Flexible boundaries in deformations of hyperbolic $3$--manifolds, Osaka J. Math. 34 (1997), 541-551.

\bibitem{fujii} M. Fujii and T. Soma, Totally geodesic boundaries are dense in the moduli space, J. Math. Soc. Japan, 49, no. 3, 1997, 589--601.

\bibitem{gard} F. P. Gardiner, Teichm\"uller theory and quadratic differentials, Pure and Applied Mathematics. A Wiley-Interscience Publication. John Wiley \& Sons, Inc., New York, 1987.

\bibitem{kapo} M. Kapovich, Hyperbolic Manifolds and Discrete Groups, Progress in Math. v. 183, Birkh\"auser, 2000.

\bibitem{kra} I. Kra, On the Nielsen--Thurston--Bers type of some self-maps of Riemann surfaces, Acta Math. 146 (1981) 231--270.

\bibitem{mask} B. Maskit, Kleinian groups, Grundlehren der mathematischen Wissenschaften 287, Springer 1988.

\bibitem{mc1} C. T. McMullen, Renormalization and $3$--manifolds which fiber over the circle, Annals of Math. Studies no. 142, Princeton, 1996.

\bibitem{morg} J. W. Morgan, On Thurston's uniformization theorem for three--dimensional manifolds, The Smith Conjecture, Academic Press, 1984, 37--125.

\bibitem{mum} D. Mumford, A remark on Mahler's compactness theorem, Proc. AMS 28 (1971) 289-294.

\bibitem{myers} R. Myers, Excellent $1$--manifolds in compact $3$--manifolds,  Topology Appl.  49  (1993),  no. 2, 115--127.

\bibitem{neureid} W. D. Neumann and A. W. Reid, Rigidity of cusps in deformations of hyperbolic $3$--orbifolds, Math. Ann. 295, 223--237 (1993).

\bibitem{stallings}  J. R. Stallings, Whitehead graphs on handlebodies. Geometric group theory down under (Canberra, 1996), 317--330,
de Gruyter, Berlin, 1999. 

\bibitem{whitehead} J. H. C. Whitehead, On certain sets of elements in a free group. Proc. London Math. Soc. 41 (1936), 48-56.

\end{thebibliography}
\end{document}